\documentclass{gENO2e}
\usepackage{ifpdf}
\usepackage{algorithm}
\usepackage{algorithmic}
\usepackage{multirow}
\usepackage{booktabs}
\usepackage{rotating}

\begin{document}


\articletype{\textbf{RESEARCH ARTICLE}}


\title{Parallel machine scheduling with step deteriorating jobs and setup times by a hybrid discrete cuckoo search algorithm}

\author{Peng Guo$^{a}$
\vspace{6pt} Wenming Cheng $^{a}$\vspace{6pt} and Yi Wang$^{b}$$^{\ast}$\thanks{$^\ast$Corresponding author. Email: ywang2@aum.edu}\\\vspace{6pt}  $^{a}${\em{School of Mechanical Engineering, Southwest Jiaotong University, Chengdu, China}};\\
$^{b}${\em{Department of Mathematics, Auburn University at Montgomery, Montgomery, AL, USA}}\\\vspace{6pt}}

\maketitle

\begin{abstract}
This article considers the parallel machine scheduling problem with step-deteriorating jobs and sequence-dependent setup times. The objective   is to minimize the total tardiness by determining the allocation and sequence of jobs on identical parallel machines. In this problem, the processing time of each job is a step function dependent upon its starting time. An individual extended time is penalized when the starting time of a job is later than a specific deterioration date. The possibility of deterioration of a job makes the parallel machine scheduling problem more challenging than ordinary ones.
A mixed integer programming model for the optimal solution is derived. Due to its NP-hard nature, a hybrid discrete cuckoo search algorithm is proposed to solve this problem. In order to generate a good initial swarm, a modified heuristic named the MBHG is incorporated into the initialization of population. Several discrete operators are proposed in the random walk of L\'{e}vy Flights and the crossover search. Moreover, a local search procedure based on variable neighborhood descent is integrated into the algorithm as a hybrid strategy in order to improve the quality of elite solutions. Computational experiments are executed on two sets of randomly generated test instances. The results show that the proposed hybrid algorithm
can yield better solutions in comparison with the commercial solver CPLEX with one hour time limit, discrete cuckoo search algorithm and the existing variable neighborhood search algorithm.

\bigskip

\bigskip
\end{abstract}


\section{Introduction}
Production scheduling as a key decision-making process has a great effect on the performance of   manufacturing and service systems. There have been a great number of research works that are involved in the classical scheduling problem. These studies assumed that the processing time of each job is constant throughout the scheduling period. However, this assumption may not be true in some real industrial settings. Examples can be found in steel production, equipment maintenance, and medical emergency, etc., where any delay or waiting in starting to process a job may incur extended time for its completion. Such kinds of jobs are called {\em deteriorating jobs}. The corresponding scheduling problem was firstly introduced by \citet{Gupta1988387} and \citet{Browne1990495}. Since then, many scheduling problems with deteriorating jobs have been extensively studied from various aspects. For details on this stream of research, \citet{Alidaee1999711}, \citet{Cheng20041} and \citet{Gawiejnowicz2008} provided comprehensive surveys of different models and problems. More recent literature that has explored scheduling problems with deteriorating jobs include \citet{Toksari2008801,Wang20091221,DSETUP_Bahalke2010749,DSETUP_Cheng20111861,Huang20111349,Bank20121,Jafari2012389,Yin2013323}
and etc.

Among various types of scheduling problems involving deteriorating jobs, those jobs with {\em step deteriorating} effect are considered in this article. A step-deteriorating job means that if a job fails to be processed prior to a pre-specified threshold, its processing time will be extended by adding an extra penalty time.
The corresponding single machine scheduling problem may date back to \citet{SD_Sundararaghavan1994394}, who firstly gave some solvable cases for minimizing the sum of the weighted completion times. Subsequently, \citet{SD_Mosheiov1995869} studied the makespan minimization problem
for both single- and multi-machine cases.
He also suggested several heuristics for all these complex problems. Subsequently, \citet{SD_Jeng2004247} considered the problem on a single machine with multiple deteriorating dates, and designed a branch and bound algorithm for deriving optimal solutions. \citet{SD_Cheng2001623} proved that the total flowtime minimization problem is NP-complete. \citet{SD_Jeng2005521} developed a branch and bound algorithm for the problem.
Owing to the intractability of the problem, \citet{SD_He20091759} {used  weighted values of the basic processing times, the deteriorating thresholds and the penalties to derive a heuristic solution.}
\citet{SD_Cheng2012928312} employed a variable neighborhood search algorithm to minimize the total completion time in a parallel machines environment. Moreover, \citet{SD_Layegh20091074} used a memetic algorithm to solve the single machine total weighted completion time problem.
 {\citet{SD_Peng2013} used a general variable neighborhood search to solve a single machine total tardiness problem.}
In particular, batch scheduling with step-deterioration has been researched \citep{SDB_Leung20081090,SDB_Mor2012587}.  In addition, a similar model which assumes that the processing time of a job is  a piecewise linear function (rather than a step function) of the waiting time is introduced by \citet{PLD_Kunnathur199056}. Other related works can be referred to \citet{PLD_Kubiak1998511,PLD_Cheng2003531,PLD_Ji200741,PLD_Wu20091230,PLD_Moslehi2010573,PLD_Farahani20131479}.

Sequence dependent setup is commonly encountered in manufacturing environments \citep{SETUP_Allahverdi1999219,SETUP_Allahverdi2008985}.
The research of  scheduling of deteriorating jobs with consideration of setup times is relatively uncommon compared with those ordinary problems with setup times.
The single machine scheduling problems with past-sequence-dependent setup times and deteriorating jobs were studied by \citet{DSETUP_Zhao2010663,DSETUP_Cheng20111861,DSETUP_Lai2011737}. In their works, the actual   processing time of a job is formulated as an   increasing function of its scheduled position and the actual processing times of jobs already processed.
In addition,  \citet{DSETUP_Bahalke2010749} suggested  a genetic algorithm and a tabu search for a single machine scheduling problem with the objective  of minimizing the makespan. Moreover, the branch and bound technique based on some dominance properties has been adopted by \citet{DSETUP_Cheng20111760,DSETUP_Lee2011782} for minimizing the maximum tardiness and the number of late jobs, respectively.
{\em To the best of our knowledge, there is no research that concerns the scheduling problem with step-deteriorating jobs and sequence dependent setup times}.

Parallel machine scheduling problem, as a generalization of the single machine problem, has been widely studied in the literature because it frequently appears in the industrial production.
In this article, we study the parallel machine scheduling problem in which the two effects of both step deterioration and sequence dependent setup times are considered.
Below the abbreviation ``PMSDST" will be used to denote the problem under consideration.
Since the minimization of the total tardiness on uniform parallel machines with sequence-dependent setup times is NP-hard \citep{SETUP_Armentano2007100}, the PMSDST, as a more general problem,  is also NP-hard.
Though the PMSDST can be formulated to a mathematical programming model, the model is impossible to be solved in a reasonable amount of computational time  when the problem size increases.
Therefore,
we propose
a hybrid meta-heuristic algorithm based on the cuckoo search to obtain a near-optimal schedule in a reasonable time.
 Computational results and comparisons to other existing algorithms demonstrate the effectiveness and the efficiency of the proposed meta-heuristic algorithm.

The rest of the article is organized as follows. The next section describes the problem in more details, and formulates a mixed integer programming model. Section \ref{sec:Sec_HCSA} proposes the search heuristic suggested in our study together with the methods to generate initial and neighborhood solutions. Then the computational results are reported and discussed in section \ref{sec:Sec_experiment}. Finally, section \ref{sec:sec_conclusion} gives some conclusions along with remarks for future research.

\section{Problem formulation\label{sec:Sec_problem}}

In the PMSDST, a set of $n$ independent jobs have to be processed on $m$ identical parallel machines. All machines are available at time zero, and each machine can handle at most one job at a time. Preemption, division, or cancelations are not allowed. For   job $j$, its actual processing time $p_{j}$ is a non-linear step function of its starting time $s_{j}$ and the deteriorating date $h_j$. If the starting time $s_j$ of job $j$ is less than or equal to the corresponding threshold $h_j$, then job $j$ only requires a basic processing time $a_j$ . Otherwise, an extra penalty $b_j$ is incurred.
Each job $j$  is given a  due dates $d_{j}$. In addition, a sequence dependent setup time for a job is also considered in this problem. That is to say, there is a setup time $\delta_{ij}$ when job $i$ precedes immediately job $j$ on the same machine.  As a matter of fact, $\delta_{jj}=0$. In addition, it is assumed that $\delta_{0j}=0$, which indicates that no setup time is needed if job $j$ is processed at the first position on a machine.
 Without loss of generality, all parameters are assumed to be integers.
The tardiness $T_{j}$ of a job $j$ can be determined by $T_{j}=\max\{0, C_{j}-d_{j}\}$, where $C_{j}$ denotes the completion time of  job $j$. Finally, the goal of the scheduling problem considered here is to find a schedule  that specifies an assignment of all jobs to various machines and the schedule of those jobs on all machines so that the total tardiness is minimized.

We next propose a mixed integer programming (MIP) model   to minimize the total tardiness for the PMSDST.
The basic structure of the MIP model is motivated by a network model formulated by \citet{MODEL_Radhakrishnan20002233}.
In order to present  constraints conveniently, two dummy jobs 0 and $n+1$ are introduced on each machine and are scheduled at the head and tail of its job sequence, and its processing time is set to 0. Additionally,
the following  notation is defined.
\begin{itemize}
  \item $u_{ij}^{k}$: the binary decision variable that is equal to 1 if job $i$ precedes immediately job $j$ on machine $k$, and 0 otherwise.

  \item $\sum_{j=1}^{n}T_j$: the total tardiness.
  \item $M$: a large positive constant such that $M\rightarrow\infty$ as $n\to \infty$. In this article, we choose  $M=\displaystyle\max_{j =1,\ldots,n}d_{j}+\displaystyle\sum_{j=1}^{n}(a_j+b_j)$.
\end{itemize}
Applying the above notations and variables, the PMSDST can be formulated as follows.

\emph{Objective function:}

\begin{equation}\label{eq:eq1}
  \mathrm{min} \quad Z=\sum_{j=1}^{n}T_{j}
\end{equation}

\emph{Subject to:}

\begin{align}
p_{j}=
  \begin{cases}
  a_{j}, & \quad s_{j}\leqslant h_j\\
  a_{j}+b_{j},& \mathrm{otherwise}
  \end{cases}, &\quad \forall j=1,\ldots,n    \label{eq:eq2}\\
\sum_{i=1}^{n}u_{0i}^{k}= 1, &\quad \forall k=1,\ldots,m \label{eq:eq3}\\
\sum_{i=1}^{n}u_{i(n+1)}^{k} = 1, &\quad \forall k=1,\ldots,m \label{eq:eq4}\\
\sum_{i=0}^{n}\sum_{k=1}^{m}u_{ij}^{k}=1, &\quad \forall j=1,\ldots,n, i\neq j \label{eq:eq5}\\
\sum_{j=1}^{n+1}\sum_{k=1}^{m}u_{ij}^{k}=1, &\quad \forall i=1,\ldots,n, i\neq j \label{eq:eq6}\\
s_{j}\geqslant \delta_{0j}+M(u_{0j}^{k}-1), &\quad \forall j=1,\ldots,n, \ \forall k=1,\ldots,m \label{eq:eq7}\\
s_{j} \geqslant C_{i}+\delta_{ij}+M (u_{ij}^{k}-1), & \quad \forall i,j=1,\ldots,n, \ \forall k=1,\ldots,m \label{eq:eq8}\\
C_{j} \geqslant s_{j}+p_{j}, & \quad \forall j=1,\ldots,n \label{eq:eq9}\\
T_{j}\geqslant C_{j}-d_{j}, &\quad \forall j=1,\ldots,n \label{eq:eq10}\\
u_{ij}^{k} \in \{0, 1\}, & \quad \forall i,j=1,\ldots,n, \ \forall k=1,\ldots,m \label{eq:eq11}\\
s_{j}, C_{j}, T_{j}\geqslant0, & \quad \forall j=1,\ldots,n \label{eq:eq12}
\end{align}

In the above  formulation, equation \eqref{eq:eq1} is the objective function that minimizes the sum of the total tardiness. Constraints \eqref{eq:eq2} define the actual processing time of a job under the consideration of step-deteriorating effect. Constraints \eqref{eq:eq3} ensure that the dummy job 0 is positioned at the beginning of the sequence before all the real jobs on each machine. Constraints \eqref{eq:eq4} guarantee that the dummy job $n+1$ is processed just after all real jobs on each machine.
Constraints \eqref{eq:eq5} and \eqref{eq:eq6} ensure the precedence relation of jobs assigned to the same machine.
Constraints \eqref{eq:eq7} state that the starting time for the first job of each machine must be equal to or greater than its initial setup time.
For other jobs $j$, the starting time $s_j$ is determined by constraints \eqref{eq:eq8} and is at least the sum of the completion time of job $i$ and the setup time from $i$ to $j$. If $i$ is not the predecessor of $j$, the subtraction of $M$ makes constraints \eqref{eq:eq8} non-restrictive.
Constraints \eqref{eq:eq9} calculate the completion time of each job. For each job $j$ ($j=1,\ldots,n$), the completion time $C_{j}$ is equal to its starting time $s_j$ plus its actual processing time $p_j$.
Constraints \eqref{eq:eq10} ensure that only the positive value of lateness can be considered as tardiness, that is, $T_j=\mathrm{max}\{0,C_j-d_j$\}. Constraints \eqref{eq:eq11} declare $u_{ij}^{k}$ as binary variables. Finally, constraints \eqref{eq:eq12} state that the starting time, completion time and tardiness of job $j$ are non-negative.


The above mathematical programming model is used to find optimal solution for our problem.
However, due to the NP-hardness of  the proposed problem, it is difficult to be solved optimally for medium- and large-sized instances. This is the main characteristic of most complex scheduling problems, which makes necessary to develop some meta-heuristic algorithms for solving the problem under study. In the next section, we propose a hybrid cuckoo search algorithm to obtain  near-optimal solutions.

\section{Hybrid discrete cuckoo search algorithm \label{sec:Sec_HCSA}}

The Cuckoo Search (CS) algorithm is one of the most new swarm intelligence optimization algorithms introduced by \citet{Yang2009210,Yang2010330}, inspired by the intelligent breeding behavior of some species of cuckoo. These species lay their eggs in the nests of host birds (of other species) with amazing abilities such as selecting the recently spawned nests and removing existing eggs so that the hatching probability of their own eggs is increased. On the other hand, some  host birds are able to engage direct contest with  infringing cuckoos. For instance, if the alien eggs are discovered by the host bird, it may either throw out these eggs or abandon the nest completely. This natural phenomenon has led to the evolution of cuckoo eggs to mimic the egg appearance of local host birds.

For simplicity in the implementation of the algorithm, the following three idealized rules are used  \citep{Yang2009210}: (1) Each cuckoo lays one egg at a time, and dumps its egg in a randomly chosen nest; (2) The best nests with high quality of eggs will carry over to the next generation; (3) The number of available host nests is fixed, and the egg laid by a cuckoo is discovered by the host bird with a probability $\rho_a\in[0, 1]$. In this case, the host bird can either discard the egg or abandon the nest so as to build a  new nest in a new location.
The last assumption can be approximated by a fraction $\rho_a$ of the $n$ nests that shall be  replaced by new ones.

The CS algorithm contains a population of $\mathcal{P}$ eggs.  {Egg $\imath \in \mathcal{P}$ is  in nest $X_{\imath}$, which corresponds to a  solution of the underlying problem. The word `nest' and the word `solution' shall be used interchangeably.
 }
In the initial phase, $\mathcal{P}$ solutions are generated randomly.
Then all of these solutions except the  one with the smallest objective value are replaced by new solutions that shall be  generated by the {\em L\'{e}vy flights}.  Specifically, if $t$ is the current iteration index, and let $\Phi^t:=\{X^t_j:j\in \mathcal{P}  \}$ be the set of all solutions at the $t$-th iteration, then
a candidate solution $\tilde{X}_\imath^{(t+1)}$ is generated by the equation
\begin{equation}\label{eq:eq_CSLF}
  \tilde{X}_\imath^{(t+1)}=X_\imath^{(t)}+\alpha_0\left(X_{\jmath}^{(t)}-X_{\imath}^{(t)}\right) \cdot t^{-\lambda},
\end{equation}
where, the constant parameter $\alpha_0$  is related to the scale of the underlying problem,  $X_{\jmath}^{(t)}$ represents a randomly selected solution at the $t$-th iteration and     $\lambda$ is a constant such that     $1<\lambda\leq3 $.

The vector $ \alpha:= \alpha_0(X_{\jmath}^{(t)}-X_{\imath}^{(t)})$ introduces a difference from the current solution $X_{\imath}^{(t)}$.
  This difference mimics the scenario in which  a cuckoo's egg is more difficult to be recognized by a host if it is  more similar to the host's egg. The term $t^{-\lambda}$ accounts for the heavy-tailed power law distribution $\text{L\'{e}vy} (\lambda)$ given by  \begin{equation}\label{eq:eq_LFPD}
  \text{L\'{e}vy} (\lambda)\sim t^{-\lambda}.
\end{equation}
The distribution $\text{L\'{e}vy} (\lambda)$
has infinite variance and mean when $1< \lambda \le 3$, thus the L\'{e}vy flight essentially provides a random walk process in hope to search for a better solution.

If the host bird identifies a cuckoo's egg in its nest (with the  probability $\rho_a$), either the alien egg shall be discarded by the host or the host shall abandon the nest. In either case, the cuckoo's egg will fail to hatch.  Thus the cuckoo has to explore a new nest (solution)  in order to successfully hatch its egg. This behavior shall be simulated by the crossover operation.  Specifically, let $\mathrm{rand}_\imath$, $\mathrm{rand}$ be two uniformly distributed random numbers on the interval $[0,1]$.  A new candidate solution  $V_{\imath}^{t+1}$ is obtained  by the equation
 \begin{equation}\label{eq:eq_CSCO}
  V_{\imath}^{t+1}=
  \begin{cases}
  \tilde{X}_{\imath}^{t+1}+\mathrm{rand} \cdot (\tilde{X}_{\gamma_1}^{t+1}-\tilde{X}_{\gamma_2}^{t+1}) & \mathrm{if}\  \mathrm{rand}_{\imath}>\rho_a\\
  \tilde{X}_{\imath}^{t+1} & \mathrm{otherwise},
  \end{cases}
\end{equation}
where, $\gamma_1$ and $\gamma_2$ are two randomly chosen indices such that $\imath$, $\gamma_1$ and $\gamma_2$ are pairwise distinct. If the candidate solution $V_{\imath}^{t+1}$ has a better quality than the old one $X_{\imath}^{t}$, that is, with a smaller objective function value,  then $V_{\imath}^{t+1}$ will be accepted as  $X_{\imath}^{t+1}$. Therefore, for $\imath\in \mathcal{P}$, we have at the $(t+1)$th iteration,
\begin{equation}\label{eqn:comparison}
X_{\imath}^{t+1}=\left\{
 \begin{array}{l l}
 V_{\imath}^{t+1}, &  \mathrm{if}\ f(V_{\imath}^{t+1})<f(X_{\imath}^{t})\\
 \tilde{X}_{\imath}^{t+1}, & \mathrm{otherwise},
 \end{array}
\right.
\end{equation}
where, $f$ is the objective function of the underlying optimization problem.
Afterwards, the best solution among all solutions found so far up to the present iteration   is recorded to memory. The same   procedure repeats until a termination criterion is met.

The CS algorithm provides more robust and precise results than   other meta-heuristic algorithms, such as  particle swarm optimization, differential evolution, or  artificial bee colony algorithms in solving the continuous optimization problems \citep{Civicioglu2013315}. Therefore, the CS algorithm  has aroused many interests and has been successfully applied to different kinds of optimization problems \citep{Walton2011710,Gandomi201317,Valian2013459,Yildiz201355}. Recently, \citet{Ouaarab2013} proposed a discrete CS algorithm to solve the traveling salesman problem. \citet{Burnwal2013951} developed a cuckoo search-based approach for scheduling optimization of a flexible manufacturing system.
\citet{Li20134732} designed a hybrid algorithm that combines the CS and a fast local search to solve a permutation flow shop scheduling problem.
%

In the remainder of this section, we propose
a hybrid discrete CS algorithm (HDCS) to solve the considered PMSDST. In order to maintain the remarkable characteristics of the standard CS, the L\'{e}vy flight equation \eqref{eq:eq_CSLF} and the crossover operator \eqref{eq:eq_CSCO} are redefined based on an integer encoding scheme. Moreover, a local search procedure based on variable neighborhood descent is employed to refine  elite solutions that are not discarded. In addition, a heuristic named the MBHG is incorporated into the random initialization to generate a population of initial solutions with certain quality.
Furthermore, to avoid the algorithm trapping into a local optimum, a restarting strategy is embedded into the search process.
The details of the proposed hybrid discrete algorithm are elaborated below.

\subsection{Solution representation}

In a meta-heuristic, the solution performance is highly dependent upon the representation of a solution.
The most commonly used solution representation for parallel machine scheduling problem is having an array of jobs for each machine that represents the processing order of  jobs assigned to that machine. Nevertheless, in this article,
we represent a solution by a    single sequence of the $n$ jobs in order to maintain the simplicity and the powerful search performance of the CS algorithm.
{\em With such a sequence,
an  actual schedule can be built simply by assigning the next unscheduled job  to the earliest available machine until all jobs are scheduled}. When multiple machines are available for an unscheduled job, the machine with the smallest index will be chosen.  The obtained sequence of jobs shall be called a solution vector or a nest vector.
Once an actual schedule is constructed, the total tardiness can be easily calculated subsequently.

\begin{figure}[htp]
  \centering
  \includegraphics[scale=1.2]{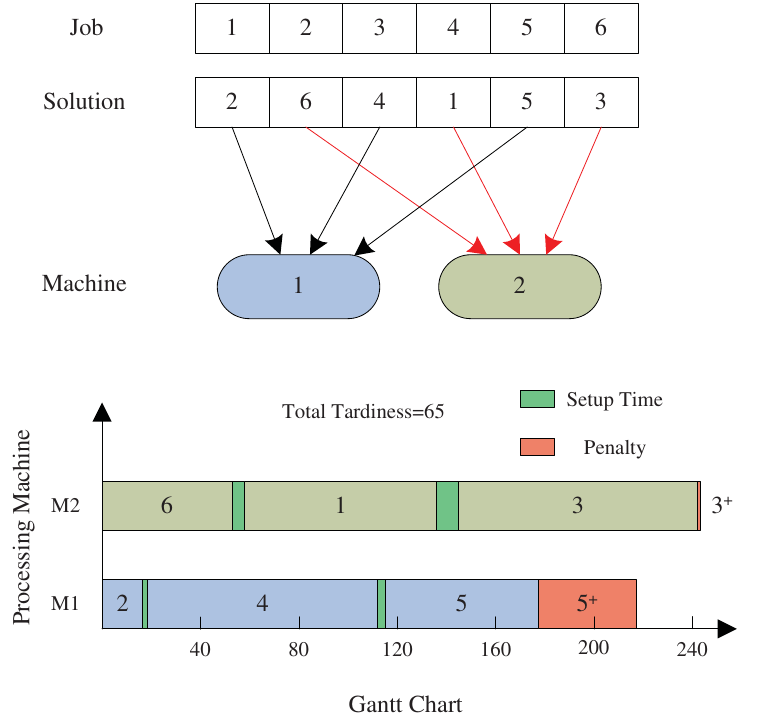}\\
  \caption{The decoding scheme of a solution}\label{Fig:fig_decoding}
\end{figure}

{To illustrate how the nest vector }is used to construct an actual schedule, below a simple example is given with {\em six} jobs and {\em two} machines. The detailed data for this problem is provided in Table \ref{tab:example_input}.
{Suppose a nest vector} for the problem  is given by $X=[2, 6, 4, 1, 5, 3]$. An actual schedule can be obtained through the aforementioned  decoding scheme, as shown in Figure  \ref{Fig:fig_decoding}. Correspondingly, the
total tardiness is 65. In  Figure \ref{Fig:fig_decoding}, jobs 2, 4 and 5 are assigned to machine 1 and jobs 6, 1 and 3 are assigned to machine 2.

\begin{table}[htp]
  \centering
  \caption{The input data of the example with six jobs and two machines}
    \begin{tabular}{lcccccc}
    \toprule
    Job($j$) & 1     & 2     & 3     & 4     & 5     & 6 \\
    \midrule
    Basic processing time $a_j$ & 78    & 17    & 97    & 93    & 62    & 53 \\
    Due date $d_j$ & 85    & 48    & 229   & 133   & 220   & 75 \\
    Deteriorating date $h_j$ & 70    & 4     & 62    & 19    & 58    & 39 \\
    Penalty $b_j$ & 18    & 33    & 1     & 17    & 40    & 31 \\
    \multirow{6}[0]{*}{Setup time $\delta_{ij}$} & 0     & 9     & 9     & 5     & 4     & 6 \\
          & 5     & 0     & 8     & 2     & 2     & 5 \\
          & 2     & 6     & 0     & 5     & 4     & 7 \\
          & 3     & 5     & 10    & 0     & 3     & 9 \\
          & 3     & 10    & 6     & 9     & 0     & 5 \\
          & 5     & 8     & 4     & 8     & 4     & 0 \\
    \bottomrule
    \end{tabular}%
  \label{tab:example_input}%
\end{table}%

\subsection{Population initialization}

In the original CS algorithm, the initial population of solutions  is often filled with $\mathcal{P}$ $n$-dimensional  vectors that are randomly generated. In order to ensure an initial population with certain quality and diversity, the BHG heuristic \citep{Biskup2008134} is modified to generate an initial solution for our problem, whereas the rest $\mathcal{P}-1$  solutions are initialized with random job permutations.

Generally, a job $j$ is likely to be scheduled in an earlier position if it has a smaller value of due date $d_j$ and a smaller value of deteriorating date $h_j$ for our problem. Therefore, in the modified BHG heuristic, a predefined list of jobs is obtained  in ascending order of their weighted values $\omega d_j+(1-\omega)h_j$, where $0<\omega<1$. Moreover, in order to reduce  computation time, we simply  assign the first $m$ jobs from the predefined list simultaneously to the $m$ machines, respectively. The detailed  modified BHG (MBHG) heuristic  algorithm is described next.
\begin{enumerate}
  \item Obtain a list of all jobs in ascending order of  the weighted value of due date and deteriorating date, i.e.  $\omega d_j+(1-\omega)h_j\leqslant  \omega d_{j+1}+(1-\omega)h_{j+1}$ for $j=1,\ldots, n-1$.
  \item The first $m$ jobs from the list are simultaneously assigned to the $m$ machines, respectively.
  \item Let $U=\{m+1, \ldots, n\}$ be the set of unscheduled jobs.
  \item Determine to which machine the job with the smallest index in $U$ should be assigned. Starting from machine 1, take the job with the smallest index in $U$ and insert it  before and after each job that is already assigned to machine 1. The insertion into the schedule of  machine 1 may start from the position after the last job and continue sequentially backward to the position before the first job assigned to machine 1.  For each instance of inserting a job into the current schedule of machine 1, calculate a new value of the total tardiness.  The same procedure is applied to machines 2, 3, \ldots, $m$ afterwards. After all the possible insertions into the schedule of each machine  are done, the algorithm chooses the sequence with the lowest total tardiness which was calculated first. Remove the first job from the set $U$.\label{step:step_BHG_insert}
  \item Repeat Step \ref{step:step_BHG_insert} until $U$ is empty.
  \item Output the final schedule on the $m$ machines and its total tardiness.
\end{enumerate}

It is found that the weight $\omega$ has a significant influence on the solution quality of the MBHG. When the weight is $\omega=0.1$, the solution provided by the MBHG is 116 for the above mentioned instance listed in table \ref{Fig:fig_decoding}. When the weight is $\omega=0.5$, the solution obtained by the MBHG is 65 for the same instance. It is hard to find the  best value of   weight $\omega$. This motivates the search with different values of $\omega$. In this study,   weight $\omega$ is varied from 0.1 to 0.9.
The best solution with the least total tardiness among the solutions obtained from different weights is selected as the final output of the MBHG heuristics.

\subsection{CS-based search}
The search process of the original CS algorithm (equations \eqref{eq:eq_CSLF}, \eqref{eq:eq_CSCO}, \eqref{eqn:comparison}) may   fall into a local optimal solution, particularly in solving a discrete optimization problem. Therefore, in order to maintain  diversity of the population, the following  search strategies are adopted.

 For the L\'{e}vy flight, a uniformly distributed random variable  $\psi$ is firstly generated in the interval (0, 1).
  Let $X^t_{\mathrm{best}}$ be the best solution found among all solutions in $\Phi^t$.
 If $\psi >0.5$, the generation of a new solutions is directed by the current best solution $X^t_{\mathrm{best}}$. Otherwise, the  searching process is around a randomly selected solution $X^t_\jmath$ {from the current population $\Phi^t$}.
  Parameter $\lambda$ is   critical   in delivering a new solution, as observed from equation \eqref{eq:eq_CSLF}.
A smaller $\lambda$ increases the global exploration ability while a large $\lambda$ tends to improve the local exploitation ability of solutions in the current search area.
Therefore, to balance the two abilities,     $\lambda$ should be linearly increased from a relatively small value to a relatively large value in the entire process of the HDCS algorithm. By doing so, the HDCS can possess good global search ability at the beginning phase of the search  while having more local search ability near the end of the iterations. In this article,  parameter $\lambda$ is chose to  vary through the  formula
\begin{equation}\label{eq:eq_callam}
  \lambda=\lambda_{\min}+\frac{\lambda_{\max}-\lambda_{\min}}{t_{\max}}t,
\end{equation}
where, $\lambda_{\min}$ and $\lambda_{\max}$ are   the initial value and the terminal value of $\lambda$ respectively, and
$t_{\max}$ is the maximum number of iterations.
In our algorithm, $\lambda_{\min}=1.1$ and $\lambda_{\max}=3$.

{Based on the above modification,  equation \eqref{eq:eq_CSLF}  of producing a new solution $\tilde{X}_{\imath}^{t+1}$ is replaced by the following equation
\begin{equation}\label{eq:eq_modLF}
  \tilde{X}_{\imath}^{t+1}=
  \begin{cases}
  X_{\imath}^{t}+\alpha_0  \cdot (X_{\imath}^{t}-X_{\mathrm{best}}^{t}) \cdot t^{-\lambda},& \psi>0.5\\
  X_{\imath}^{t}+\alpha_0  \cdot (X_{\imath}^{t}-X_{\jmath}^{t}) \cdot t^{-\lambda}, & \mathrm{otherwise}.
  \end{cases}
\end{equation}
In this article, we set $\alpha_0=1$.

Since the standard CS algorithm is originally designed for solving a continuous optimization problem, the operations such as subtraction, multiplication and addition involved in the L\'{e}vy flight  must be redefined for a discrete scheduling problem. These operations shall be  redefined by  some neighborhood search strategies to suite a discrete problem.
Next, the discrete operators used in equation \eqref{eq:eq_modLF} are described in detail.
\begin{itemize}
  \item Subtraction
\end{itemize}

  Our algorithm involves one subtraction operation between two solutions, say,  the minuend $X_1$ and the
  subtrahend  $X_2$. Since the solution is a job permutation, a comparison strategy is carried out for producing a new solution chain $X'$. If the given elements of two solutions are identical, the corresponding element of the new chain is set equal to the number zero. Otherwise,
  the element of the new chain is copied from the corresponding one of $X_1$.
  Figure \ref{fig:fig_subtract} shows the redefined subtraction operation.

  \begin{figure}[htp]
  \centering
  \includegraphics[scale=1]{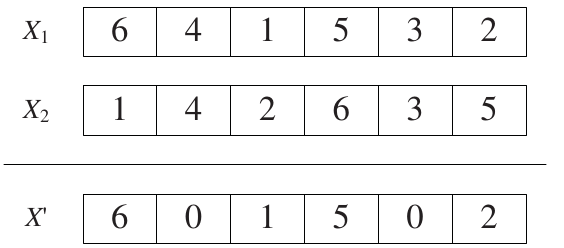}\\
  \caption{The subtraction operation $X_1-X_2$.}\label{fig:fig_subtract}
\end{figure}

%
\begin{itemize}
  \item Multiplication
\end{itemize}

  The multiplication operator carries out the product of  the scalar $\sigma:=t^{-\lambda}$ and  a solution chain $X'$ obtained by the subtraction operation.
   Let $X''$ denote the  output of the multiplication.
   In the redefined manner of the multiplication operation,
  a  uniformly distributed random number between 0 and 1 is generated for each element in the chain $X'$.
 %
  The random number is then compared with  $\sigma$.
  If the random number is more than or equal to $\sigma$, the element in the offspring chain $X''$ is set equal to the element of the same position in the chain $X'$; otherwise, its corresponding position  is occupied by the number zero.
%
%
  Figure \ref{fig:fig_multiply} represents the process of the multiplication operation on the scalar $\sigma$  and a solution chain.

  \begin{figure}[htp]
  \centering
  \includegraphics[scale=1]{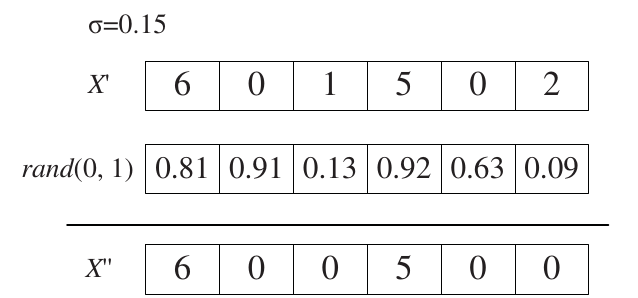}\\
  \caption{The multiplication operation}\label{fig:fig_multiply}
\end{figure}
\begin{itemize}
  \item Addition
\end{itemize}

  The addition operation is modified by a number of  swap moves. Two objects manipulated in the operation are the current solution $X$ and the chain $X''$ produced from the multiplication. Once the addition operation is finished, a new  solution $X_{\mathrm{new}}$ is generated. Initially, $X_{\mathrm{new}}$ is set  equal to the incumbent  $X$. Then
  for each element  $X''(i)$ of $X''$, if $X''(i)$ is more than 0,
{ the job  in position $X''(i)$  of $X_{\mathrm{new}}$ and  the job in position $X_{\mathrm{new}}(i)$ of $X_{\mathrm{new}}$ are swapped.}
Figure \ref{fig:fig_add} illustrates the manner in which the addition performs.

\begin{figure}[htp]
  \centering
  \includegraphics[scale=1]{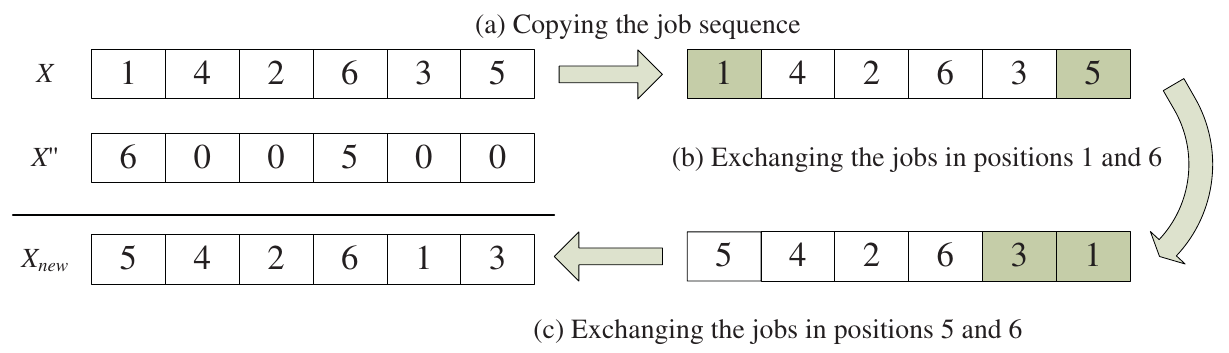}\\
  \caption{The addition operation}\label{fig:fig_add}
\end{figure}

As it can be seen, the L\'{e}vy flight is redefined by the above three discrete operators. Likewise, the crossover operation described by equation \eqref{eq:eq_CSCO} is replaced by the {\em order crossover operator}.
The {order crossover} \citep{Davis1985162_GAOC} has shown     robust performance in the genetic algorithm.
For each solution,
a random  number is generated and compared with the probability $\rho_a$ to determine whether the crossover operator should be executed. If the random number is greater than $\rho_a$, the solution
 is selected as a parent for producing a new offspring. Another solution
  is randomly chosen from the population of solutions and regarded as a second parent. With the two parent solutions chosen, the order crossover operation is performed as below.
\begin{enumerate}
  \item Select randomly  a subsequence  from one parent, called the first parent $X_1$.
  \item Produce a proto-child by copying the subsequence  into its corresponding position.
  \item Check all the jobs of the other parent, called the second parent $X_2$, from left to right. If a job is already assigned in the subsequence from $X_1$, then skip to the next job in $X_2$ that is not scheduled yet and assign it to the first available position of the new solution. The jobs taken from $X_2$ are the ones that the proto-child needs.
\end{enumerate}
The above crossover operator will produce two offspring solutions because each parent may be chosen as the first parent.  The incumbent solution shall be updated by the better one of the two offsprings {in the next generation.}
The procedure is illustrated in figure \ref{fig:fig_crossover}.

\begin{figure}[htp]
  \centering
  \includegraphics[scale=1]{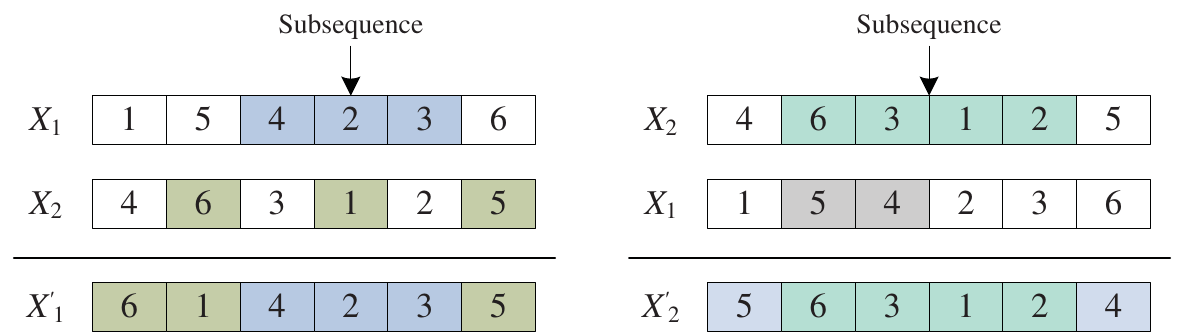}\\
  \caption{The order crossover operation}\label{fig:fig_crossover}
\end{figure}

\subsection{Local search}

To improve the performance,
the HDCS implements variable neighborhood descent (VND) as a local search strategy for a set of elite solutions.
The elite set is constructed by finding the best $\tau$ number of solutions from the population in each iteration.  The preset value of  $\tau$ can be chosen by the following equation
\begin{equation}\label{eq:eq_tau}
  \tau=\max\{3, \mathrm{round}[\mathcal{P}\times (1-\rho_a)]\}.
\end{equation}
{Equation \eqref{eq:eq_tau} ensures that there are  at least 3 elite solutions. If there are too few   elite solutions, the local search strategy can not improve the performance of the algorithm   effectively. On the other hand,   too many  elite solutions may consume too much computational time.}
Each   elite solution
indicates a promising region  to search for the optimal solution.
  Hence, the VND local search increases the possibility of obtaining an optimal solution from diverse areas of the solution space.
In our implementation, we use three neighborhood structures, which are based on \emph{swap}, \emph{insert} and \emph{inverse} operations.
The  {swap} operation is done by randomly selecting the $i$th and $j$th jobs from an elite solution and then swapping them directly. The  {insert} operation is done by choosing   job $i$ at random from a solution and inserting it into the position immediately preceding a randomly chosen  job $j$ of the solution. The {inverse} operation is done by inverting the subsection between two  randomly chosen jobs of a solution. {Each of the  three neighborhood operations is repeated $n$ times when it is executed. An example of each operation is shown in Figure} \ref{fig:fig_vns}.
\begin{figure}[htp]
  \centering
  \includegraphics[scale=1]{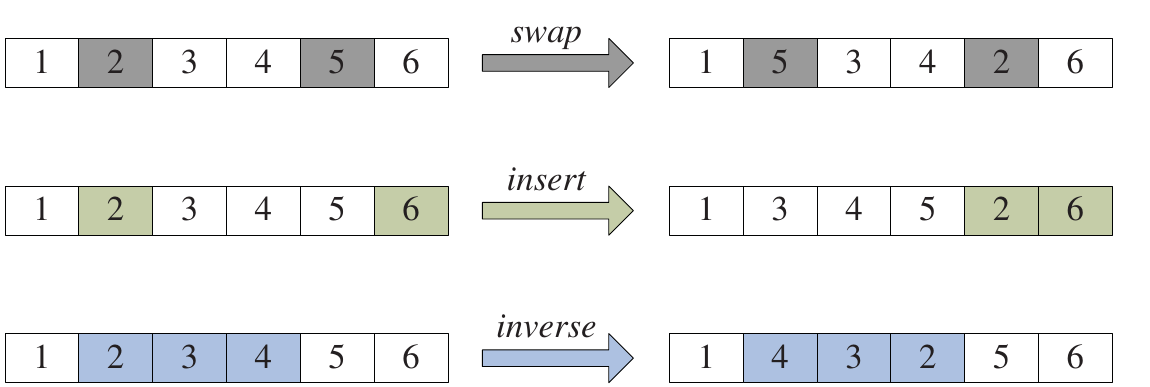}\\
  \caption{The \emph{swap}, \emph{insert} and \emph{inverse} operations for local search}\label{fig:fig_vns}
\end{figure}

The proposed variable neighborhood descent is applied to enhance the quality of an elite solution.
Each  elite solution is regarded as a seed permutation.
The search process repeats each neighborhood for the best local optimum until no improvement appears.
If the local optimum is better than the seed solution, update the seed solution and continue the search in the neighborhood; otherwise move to the next neighborhood and continue to search there.
The pseudo-code of the suggested local search is summarized in Algorithm \ref{alg:alg_LS}.

\begin{algorithm}[htp]
\caption{\label{alg:alg_LS} Local search}
\begin{algorithmic}[1]
\REQUIRE Elite solutions for the HDCS
\ENSURE Improved elite solutions for the HDCS
\FOR {each solution $X$ in the elite set}
\STATE $\kappa$=1;
\WHILE {$\kappa \leqslant 3$}
\IF {$\kappa==1$}
\STATE Apply \emph{Swap} operations $n$ times to $X$  to obtain a local optimum $X^{'}$;
\ELSIF {$\kappa==2$}
\STATE Apply \emph{Insert} operations $n$ times to $X$  to obtain a local optimum $X^{'}$;
\ELSIF {$\kappa==3$}
\STATE Apply \emph{Inverse} operations $n$ times to $X$  to obtain a local optimum $X^{'}$;
\ENDIF
\IF {$X^{'}$ is better than $X$}
\STATE $X\leftarrow X^{'}$;
\STATE Continue the $k$th neighborhood search operation;
\ELSE
\STATE $\kappa=\kappa+1$;
\ENDIF
\ENDWHILE
\ENDFOR
\end{algorithmic}
\end{algorithm}


\subsection{Restarting strategy and stopping criterion}
To maintain the diversity of the population of solutions, it is common to use restarting strategy in evolutionary algorithms \citep{Ruiz2006461}. In our proposed algorithm,  partial solutions adopt a restarting mechanism  at each iteration. Firstly, the population is sorted in ascending order of objective function values. Then the first 90\% solutions from the sorted list are retained in the population. The remaining 10\%  solutions shall be  replaced with newly generated nests  {by random job permutations}.

Regarding the stopping criterion, two commonly used termination conditions are employed in our proposed algorithm. That is to say, if the number of iterations $t$ exceeds the maximum number of iterations $t_{\max}$ or the best solution $X_{\mathrm{best}}$ is not improved upon in $t_{\mathrm{nip}}$ successive iterations, the procedure is terminated; otherwise, a new iteration will start and the iteration counter $t$ is increased by one. When the algorithm terminates, the best found  solution $X_{\mathrm{best}}$ and its total tardiness are output.

\subsection{Outline of the HDCS}

The search framework of the proposed HDCS algorithm for solving the PMSDST is shown in Figure \ref{Fig:HCS}. The algorithm starts with an initial population of $\mathcal{P}$ solutions, among which one solution is generated by the MBHG, and   the rest is generated by random permutations of $n$ jobs. Then the solution $X_{\mathrm{best}}$ with the smallest objective value is identified by calculating the function value of each solution.
Subsequently, the algorithm iterates search for a better solution until a stopping condition is met. At each iteration, the first $\tau$ elite solutions are determined according to equation \eqref{eq:eq_tau} by their smaller objective function values. The remaining solutions are regarded as  {\em normal solutions}.
For each of the normal solutions, the modified L\'{e}vy flight by equation \eqref{eq:eq_modLF} is performed to search for a new  solution either around   $X_{\mathrm{best}}$ achieved in the previous iteration or a randomly selected neighbor.
The elite solutions are improved by a local search based on the variable neighborhood descent.
Then, a fraction $\rho_a $ of the   solutions are replaced by the  new candidate solutions   produced by performing the order crossover operation. To avoid  premature convergence, the restarting strategy is utilized for the  10\% worst solutions. All the solutions are evaluated for their objective values and the $X_{\mathrm{best}}$ is updated if a new better solution is found. The algorithm is then continued to next iteration.  When the algorithm is terminated, the  best solution found by the algorithm and its objective value are output for the problem under study.


\begin{figure}[htp]
  \centering
  \includegraphics[scale=0.6]{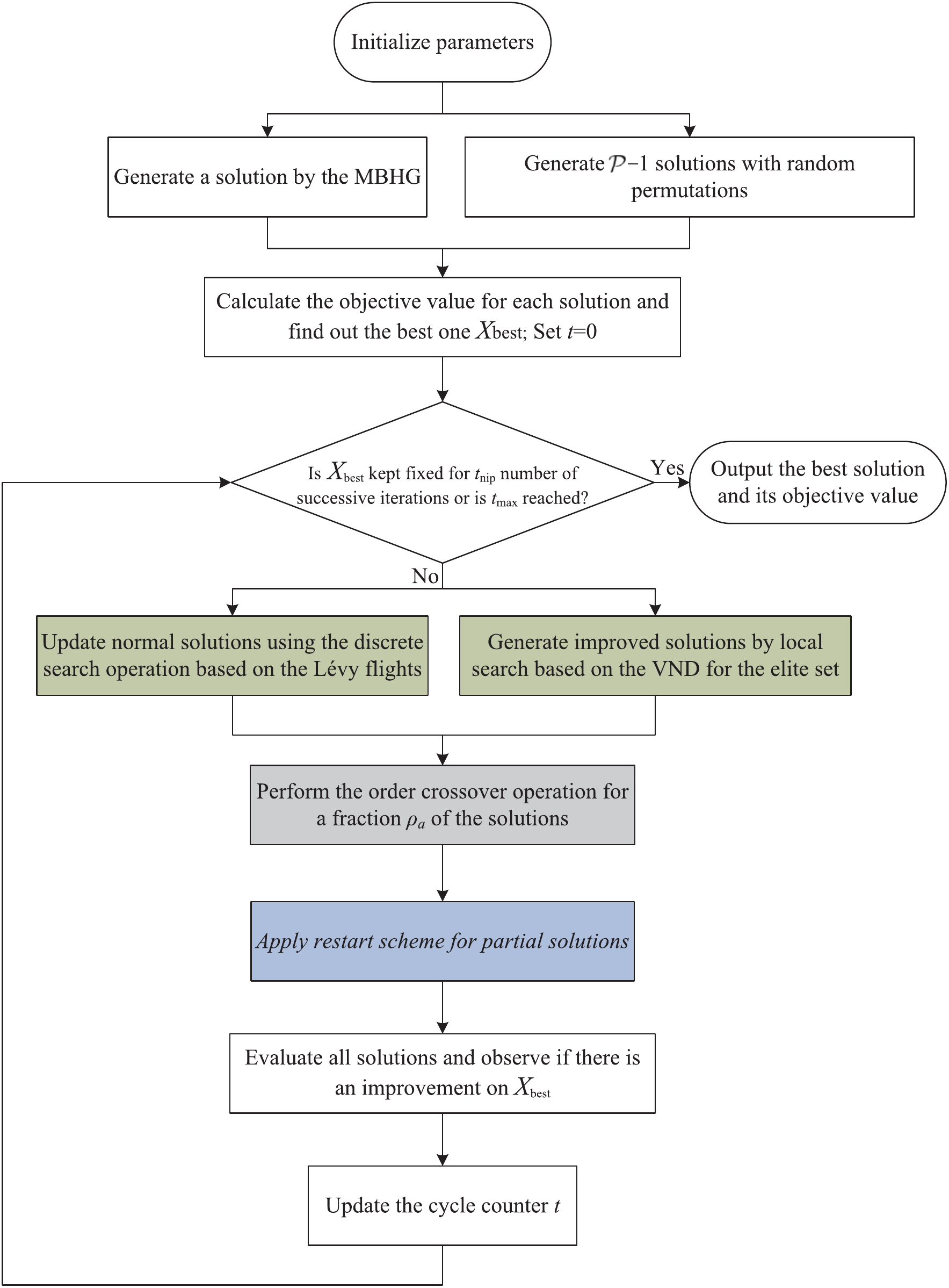}\\
  \caption{The flowchart for the proposed hybrid cuckoo search algorithm}\label{Fig:HCS}
\end{figure}


\section{Computational experiments\label{sec:Sec_experiment}}
There is no study on parallel machine scheduling problem with step-deterioration jobs and setup times in the existing literature. Therefore, our article seems to be the first study of this problem.
In order to evaluate the effectiveness of the proposed algorithm, the performance of the HDCS is compared with the discrete cuckoo search algorithm without local search strategy (DCS) and the variable neighborhood search proposed by \citet{SD_Cheng2012928312} (VNS).
The three algorithms are coded in Matlab 7.11 and executed on an Intel Pentium dual-core 2.6 GHz PC with 2GB RAM under Windows XP environment.
The performance of each algorithm is measured by
the relative percentage deviation ($RPD$) calculated by the formula
\begin{equation}\label{eq:eq_RPD}
  RPD(\%)=\frac{Alg_{sol}-Min_{sol}}{Min_{sol}}\times 100
\end{equation}
for each test instance,
where $Alg_{sol}$ is the objective function value obtained for a given algorithm and
 $Min_{sol}$ is the best solution of all experiments for each  test problem.

\subsection{Data generation of test problems\label{subsec:subsec_dategeneration}}

In order to analyze the performance of our algorithm for the scheduling problem under consideration, different sizes of test instants are needed.
In this section,
two sets of instances are defined.
For   small-sized instances, the following combinations of number $n$ of jobs and number $m$ of machines are tested: $n=\{8, 10, 12, 14\}$ and  $m=\{2, 3, 4\}$. For   large instances, the combinations are $n=\{30, 40, 50, 60\}$ and $m=\{4, 6, 8\}$.
 In all cases,   the following rules are used to generate  job data.

 The basic processing times $a_j$ of job $j$ is generated from a discrete  uniform distribution on the  interval [1, 100]. The penalty $b_j$ is assumed to be drawn from the uniform distribution $U$[1, 100$\times \xi$], where $\xi=0.5$ is chosen in our study. Set $\beta=\sum_{j=1}^{j=n}{a_j}/m$. The deteriorating dates are generated randomly by a discrete uniform distribution on  each of the intervals $\mathrm{H}_1:=[1,  0.5\times\beta]$, $\mathrm{H}_2:=[0.5\times\beta, \beta]$ and $\mathrm{H}_3:=[1, \beta]$. The setup times are produced according to a uniform distribution in the range [1, 10]. Moreover, the due dates of jobs are randomly generated from a uniform distribution on the interval [1, $\bar{C}_{\max}$], where $\bar{C}_{\max}$ is the value of the maximal completion time obtained by the following rule:
first of all, a job permutation is achieved by arranging the jobs in the non-decreasing order of the ratio $a_j/b_j$. Then, a corresponding schedule with no idle times can be constructed by assigning  unscheduled jobs  to the earliest available machine. Once a schedule is built, the maximal completion time $\bar{C}_{\max}$  is calculated.
{For each possible combination,}
 10 random replicates are generated. Therefore,
there are 360 small instances and 360 large instances to test the algorithms.

\subsection{Parameter calibration}

Parameter selection can significantly affect the quality of a algorithm. In this study, suitable parameter values were selected taking into consideration both solution quality and computational efficiency. To determine   appropriate values of   parameters,   extensive computational tests were performed on 27  instances, including both small-sized and large-size ones.   {These instances are generated from the above defined rule} with the following parameter values: $\mathcal{P}$=15,20,25,30,35,40; $\rho_a$=0.6,0.7,0.8,0.9; $t_{\max}$=100,200,300,400; $t_{\mathrm{nip}}$=30,50,70. The preliminary tests using these parameter values  show that  the following set of parameter values seems to provide the best performance within a reasonable computational time: $\mathcal{P}=30$, $\rho_a=0.8$, $t_{\max}=200$ and $t_{\mathrm{nip}}=50$. Hence, this set of parameter values will be adopted for all subsequent experiments in this study.

\subsection{Results and discussion}
In this section, we perform  comparative evaluation of the proposed HDCS, the VNS, the DCS, and a commercial solver CPLEX 12.5 based on the branch-and-cut algorithm on the same benchmark of instances explained in Section \ref{subsec:subsec_dategeneration}. Due to the intractability of the studied problem, the CPLEX 12.5 is only used to solve the MIP model with small instances. A 3600-second time limit is imposed.  When solved by the CPLEX, a particular run is simply terminated and the best current integer solution is returned if the optimal solution has not been found or verified in that amount of time.

Owing to the stochastic nature of the evaluated meta-heuristics, each algorithm is run {\em  five times} for a given instance to reach  reliable results. 
The computational results of the best ($\mathrm{RPD}_B$), average ($\mathrm{RPD}_A$) and worst ($\mathrm{RPD}_W$) deviation of each combination of $n$ (number of jobs) and $m$ (number of machines) obtained by the VNS, the DCS, the HDCS are summarized  in tables \ref{tab:tab_ARPD_SI} and \ref{tab:tab_ARPDRT_LI}. Recall that  for each combination of $n$ and $m$ in our experiment, 10 random replicates are generated and for each instance the algorithm is run 5 times. Thus corresponding to each combination of $n$ and $m$, each  value of $\mathrm{RPD}_B$ and $\mathrm{RPD}_W$ is found in the following way: firstly the best (respectively, the worst) RPD value is found among the five runs of an algorithm applied to a particular instance, then the average value of the 10 best (respectively, the worst) RPD values of those 10 random instances is recorded as $\mathrm{RPD}_B$ (respectively, $\mathrm{RPD}_W$) in tables \ref{tab:tab_ARPD_SI} and \ref{tab:tab_ARPDRT_LI}. Likewise, the value of $\mathrm{RPD}_A$ is obtained in the following way: firstly the average RPD value is found for the five runs of an algorithm applied to a particular instance, then the average value of those 10 average  RPD values of the 10 random instances is recorded as $\mathrm{RPD}_A$ in tables \ref{tab:tab_ARPD_SI} and \ref{tab:tab_ARPDRT_LI}.

Table \ref{tab:tab_ARPD_SI} reports the experimental results of the small instances. Column 1 shows the combination of $n$ and $m$ for an instance. Columns 2--5 report the computational results of the four methods, respectively. If the  CPLEX can yield the optimal solution, the RPD value is calculated over the optimal total tardiness.
 {The CPLEX solver is able to deliver    optimal solutions for all instances with eight jobs.
For the 10-job instances, the CPLEX is able to optimally solve all the instances with three and four machines, but only  23 of 30 instances with two machines. }
Regarding the 12- and 14-job cases, the CPLEX only gives   optimal solutions for a small portion of instances. In total, the CPLEX can yield   optimal solutions for 223 of the 360 small instances.

The results revealed by table  \ref{tab:tab_ARPD_SI} is exciting. We see that the HDCS can deliver the best
 {  $\mathrm{RPD}_B$ values except for the instances $10\times 2$ with $\mathrm{H}_1$ and  $14\times 3$ with $\mathrm{H}_2$.} Even for these two instances, the
 {  $\mathrm{RPD}_B$ values given by the HDCS are only slightly larger than the corresponding best values given by one of other three algorithms.
The VNS is efficient in solving a parallel machine flow time problem with step-deterioration \citep{SD_Cheng2012928312}, however, it does not perform as well as the HDCS  for the PMSDST problem. This is because the consideration of setup times and the total tardiness criterion in the PMSCST problem makes the VNS   more easily get trapped into a  local optimal solution.

In order to  further analyze   the results, the one-way analysis of variance (ANOVA) is adopted to check whether the observed difference in the $\mathrm{RPD}_{A}$ values for different algorithms are statistically significant when  one of the  deteriorating intervals ($\mathrm{H}_1$, $\mathrm{H}_2$ and $\mathrm{H}_3$) is applied. The means plots  along with the Tukey Honestly Significant Difference (HSD) intervals at the 95\% confidence level are given in figure \ref{fig:fig_RPDSI_ANOVA} for each algorithm applied to each of the  three different deteriorating intervals. The means obtained by the four algorithms for the  same deteriorating interval (one of $\mathrm{H}_1$, $\mathrm{H}_2$ and $\mathrm{H}_3$) are connected by a line segment.
{If the Tukey HSD intervals of two algorithms for the same deteriorating interval have overlapping, then the performances of the two  algorithms are not statistically significantly different.}
It is clearly seen from figure \ref{fig:fig_RPDSI_ANOVA}, that the proposed HDCS is statistically better than  other methods by having its Tukey HSD intervals of the RPD values lower than those of other methods when applied to  each of the three intervals $\mathrm{H}_1$, $\mathrm{H}_2$ and $\mathrm{H}_3$.

The computation times of the four methods   for each combination of $n$ (number of jobs) and $m$ (number of machines) are also presented in table \ref{tab:tab_RT_SI}. They are obtained in a similar manner as to obtain the $\mathrm{RPD}_A$. When the size of instance increases, the computation time of the CPLEX grows significantly. Once the number of jobs exceeds 10, the CPLEX completely exhausts the given time limit for most instances. The   DCS consumes the least time comparing with other algorithms, but its performance is the worst  among all algorithms. Alternatively, the computation times of the VNS and the HDCS are very modest.
\begin{table}[htp]\footnotesize
  \centering
  \caption{ Average Relative Percentage Deviation (RPD) of small instances for the algorithms}
    \begin{tabular}{cccccccccccc}
    \toprule
    \multicolumn{2}{c}{\multirow{2}[4]{*}{Instance}} & \multirow{2}[4]{*}{MIP} & \multicolumn{3}{c}{VNS}      & \multicolumn{3}{c}{DCS} & \multicolumn{3}{c}{HDCS} \\
    \cmidrule[0.05em](r){4-6}
    \cmidrule[0.05em](lr){7-9}
    \cmidrule[0.05em](lr){10-12}
    \multicolumn{2}{c}{} &       & $\mathrm{RPD}_B$  & $\mathrm{RPD}_A$  & $\mathrm{RPD}_W$  & $\mathrm{RPD}_B$  & $\mathrm{RPD}_A$  & $\mathrm{RPD}_W$  & $\mathrm{RPD}_B$  & $\mathrm{RPD}_A$  & $\mathrm{RPD}_W$ \\
    \midrule
    \multirow{12}[2]{*}{$\mathrm{H}_1$} & 8$\times$2   & 0.00  & 0.00  & 0.22  & 0.31  & 0.31  & 1.63  & 4.79  & 0.00  & 0.00  & 0.00  \\
          & 8$\times$3   & 0.00  & 0.00  & 0.00  & 0.00  & 0.00  & 1.50  & 5.72  & 0.00  & 0.00  & 0.00  \\
          & 8$\times$4   & 0.00  & 0.00  & 0.00  & 0.00  & 0.00  & 0.72  & 1.21  & 0.00  & 0.00  & 0.00  \\
          & 10$\times$2  & 0.00  & 0.44  & 0.96  & 1.99  & 6.43  & 11.13  & 14.54  & 0.05  & 0.05  & 0.05  \\
          & 10$\times$3  & 0.00  & 0.00  & 0.13  & 0.45  & 1.14  & 5.41  & 7.35  & 0.00  & 0.01  & 0.08  \\
          & 10$\times$4  & 0.00  & 0.00  & 0.34  & 1.50  & 0.58  & 3.83  & 5.79  & 0.00  & 0.00  & 0.00  \\
          & 12$\times$2  & 0.36  & 0.00  & 0.72  & 2.32  & 15.05  & 18.06  & 19.87  & 0.00  & 0.00  & 0.00  \\
          & 12$\times$3  & 0.00  & 0.14  & 1.76  & 5.83  & 14.26  & 17.11  & 18.97  & 0.00  & 0.09  & 0.22  \\
          & 12$\times$4  & 0.00  & 0.24  & 2.58  & 5.18  & 56.32  & 58.86  & 60.87  & 0.00  & 0.02  & 0.03  \\
          & 14$\times$2  & 12.34  & 0.00  & 2.74  & 8.02  & 53.13  & 55.09  & 55.80  & 0.00  & 0.41  & 0.92  \\
          & 14$\times$3  & 1.13  & 0.12  & 1.44  & 3.13  & 19.91  & 23.26  & 24.60  & 0.00  & 0.26  & 0.57  \\
          & 14$\times$4  & 1.51  & 0.02  & 1.14  & 2.71  & 19.60  & 20.75  & 21.83  & 0.00  & 0.02  & 0.11  \\
            \multicolumn{2}{c}{Average} & 1.28  & 0.08  & 1.00  & 2.62  & 15.56  & 18.11  & 20.11  & 0.00  & 0.07  & 0.17  \\
            &       &       &       &       &       &       &       &       &       &       &  \\
    \multirow{12}[2]{*}{$\mathrm{H}_2$} & 8$\times$2  & 0.00  & 0.00  & 0.29  & 0.59  & 0.11  & 0.58  & 0.75  & 0.00  & 0.00  & 0.00  \\
          & 8$\times$3  & 0.00  & 0.00  & 0.00  & 0.00  & 0.00  & 0.04  & 0.18  & 0.00  & 0.00  & 0.00  \\
          & 8$\times$4  & 0.00  & 0.00  & 0.00  & 0.00  & 0.00  & 0.00  & 0.00  & 0.00  & 0.00  & 0.00  \\
          & 10$\times$2  & 0.00  & 0.00  & 0.23  & 0.65  & 1.39  & 1.64  & 1.99  & 0.00  & 0.00  & 0.00  \\
          & 10$\times$3  & 0.00  & 0.00  & 0.15  & 0.33  & 7.00  & 8.24  & 9.96  & 0.00  & 0.00  & 0.00  \\
          & 10$\times$4  & 0.00  & 0.00  & 0.04  & 0.06  & 1.03  & 4.51  & 6.83  & 0.00  & 0.00  & 0.00  \\
          & 12$\times$2  & 0.87  & 0.00  & 2.73  & 10.07  & 25.28  & 26.65  & 27.42  & 0.00  & 0.09  & 0.23  \\
          & 12$\times$3  & 0.23  & 0.18  & 1.97  & 4.18  & 8.00  & 9.11  & 9.66  & 0.00  & 0.39  & 1.05  \\
          & 12$\times$4  & 0.00  & 0.00  & 1.02  & 1.86  & 4.87  & 6.94  & 8.76  & 0.00  & 0.10  & 0.33  \\
          & 14$\times$2  & 4.58  & 0.03  & 1.15  & 3.13  & 31.56  & 36.69  & 38.33  & 0.00  & 0.58  & 2.29  \\
          & 14$\times$3  & 0.28  & 0.00  & 1.17  & 3.13  & 30.89  & 36.37  & 38.96  & 0.09  & 0.49  & 1.63  \\
          & 14$\times$4  & 0.33  & 0.85  & 2.31  & 4.69  & 16.57  & 17.47  & 18.01  & 0.00  & 0.96  & 2.08  \\
           \multicolumn{2}{c}{Average} & 0.52  & 0.09  & 0.92  & 2.39  & 10.56  & 12.35  & 13.40  & 0.01  & 0.22  & 0.63  \\
           &       &       &       &       &       &       &       &       &       &       &  \\
          & 8$\times$2  & 0.00  & 0.00  & 0.00  & 0.00  & 0.00  & 0.97  & 1.70  & 0.00  & 0.00  & 0.00  \\
          & 8$\times$3  & 0.00  & 0.00  & 0.00  & 0.00  & 0.00  & 0.37  & 0.93  & 0.00  & 0.00  & 0.00  \\
          & 8$\times$4  & 0.00  & 0.00  & 0.00  & 0.00  & 0.00  & 0.07  & 0.28  & 0.00  & 0.00  & 0.00  \\
    \multirow{9}[0]{*}{$\mathrm{H}_3$} & 10$\times$2  & 0.00  & 0.00  & 0.42  & 0.79  & 7.86  & 10.43  & 11.33  & 0.00  & 0.00  & 0.00  \\
          & 10$\times$3  & 0.00  & 0.00  & 0.82  & 2.07  & 3.93  & 5.92  & 9.85  & 0.00  & 0.00  & 0.00  \\
          & 10$\times$4  & 0.00  & 0.00  & 0.08  & 0.25  & 1.54  & 2.08  & 2.77  & 0.00  & 0.00  & 0.00  \\
          & 12$\times$2  & 1.00  & 0.35  & 3.78  & 5.76  & 54.08  & 61.27  & 65.25  & 0.00  & 0.04  & 0.14  \\
          & 12$\times$3  & 0.20  & 0.20  & 3.61  & 8.81  & 21.39  & 26.33  & 29.44  & 0.00  & 0.00  & 0.00  \\
          & 12$\times$4  & 0.11  & 0.04  & 1.18  & 3.09  & 6.84  & 10.42  & 12.23  & 0.04  & 0.35  & 1.06  \\
          & 14$\times$2  & 0.61  & 2.69  & 5.18  & 8.25  & 38.65  & 46.73  & 52.13  & 0.00  & 1.48  & 3.14  \\
          & 14$\times$3  & 1.38  & 0.00  & 2.70  & 5.56  & 26.37  & 34.83  & 39.53  & 0.00  & 0.10  & 0.52  \\
          & 14$\times$4  & 0.29  & 0.67  & 1.98  & 4.04  & 18.38  & 23.82  & 27.01  & 0.00  & 0.60  & 0.95  \\
          \multicolumn{2}{c}{Average} & 0.30  & 0.33  & 1.65  & 3.22  & 14.92  & 18.60  & 21.04  & 0.00  & 0.21  & 0.48  \\
    \bottomrule
    \end{tabular}%
  \label{tab:tab_ARPD_SI}%
\end{table}%

\begin{table}[htp]\footnotesize
  \centering
  \caption{Computation times of small instances for the algorithms ( in seconds)}
    \begin{tabular}{ccccccccccccc}
    \toprule
    \multirow{2}[2]{*}{Instance} & \multicolumn{4}{c}{$\mathrm{H}_1$} & \multicolumn{4}{c}{$\mathrm{H}_2$}        & \multicolumn{4}{c}{$\mathrm{H}_3$} \\
      \cmidrule[0.05em](r){2-5}
    \cmidrule[0.05em](lr){6-9}
    \cmidrule[0.05em](lr){10-13}
          & MIP   & VNS   & DCS   & HDCS  & MIP   & VNS   & DCS   & HDCS  & MIP   & VNS   & DCS   & HDCS \\
          \midrule
    8$\times$2   & 5.22  & 1.06  & 0.52  & 1.56  & 30.91  & 1.08  & 0.55  & 1.56  & 16.12  & 1.09  & 0.63  & 1.55  \\
    8$\times$3   & 2.01  & 0.99  & 0.58  & 1.50  & 1.17  & 1.00  & 0.51  & 1.49  & 2.83  & 0.99  & 0.56  & 1.46  \\
    8$\times$4   & 0.33  & 0.96  & 0.49  & 1.38  & 0.20  & 0.93  & 0.49  & 1.37  & 0.31  & 0.95  & 0.41  & 1.38  \\
    10$\times$2  & $1188.28^1$  & 1.94  & 0.83  & 2.62  & $2140.39^4$  & 1.95  & 0.62  & 2.60  & $1590.56^2$  & 1.97  & 0.67  & 2.72  \\
    10$\times$3  & 209.54  & 1.80  & 0.81  & 2.60  & 423.81  & 1.77  & 0.64  & 2.46  & 97.56  & 1.75  & 0.74  & 2.44  \\
    10$\times$4  & 14.91  & 1.70  & 0.78  & 2.28  & 6.29  & 1.62  & 0.59  & 2.25  & 48.63  & 1.70  & 0.71  & 2.32  \\
    12$\times$2  & $3559.68^9$  & 3.00  & 1.01  & 4.34  & $3259.20^9$  & 2.80  & 0.70  & 4.15  & $3415.89^9$  & 3.02  & 0.80  & 4.58  \\
    12$\times$3  & $2496.76^6$  & 2.87  & 0.95  & 4.02  & $3101.88^8$  & 2.93  & 0.75  & 4.20  & $2328.34^5$  & 2.86  & 0.83  & 4.00  \\
    12$\times$4  & $1873.83^5$  & 2.78  & 0.89  & 3.91  & $1160.38^3$  & 2.76  & 0.72  & 3.61  & $1826.88^4$  & 2.69  & 0.80  & 3.80  \\
    14$\times$2  & $3600.00^{10}$  & 4.73  & 1.01  & 6.74  & $3600.00^{10}$  & 4.53  & 0.90  & 6.82  & $3258.91^9$  & 4.81  & 1.00  & 6.59  \\
    14$\times$3  & $3600.00^{10}$  & 4.59  & 0.85  & 6.88  & $2791.57^7$  & 4.62  & 0.82  & 6.61  & $2862.13^7$  & 4.41  & 0.91  & 6.28  \\
    14$\times$4  & $2576.57^6$  & 4.10  & 0.96  & 5.69  & $2755.64^7$  & 4.24  & 0.95  & 6.24  & $3039.72^7$  & 4.26  & 0.86  & 6.23  \\
          &       &       &       &       &       &       &       &       &       &       &       &  \\
    Average & 1593.93  & 2.54  & 0.81  & 3.63  & 1605.95  & 2.52  & 0.69  & 3.61  & 1540.66  & 2.54  & 0.74  & 3.61  \\
    \bottomrule
    \multicolumn{13}{c}{$*^\ell$ indicates that $\ell$ out of 10 instances can not be solved optimally by the CPLEX within the 3600 seconds time limit.}
    \end{tabular}%
  \label{tab:tab_RT_SI}%
\end{table}%

\begin{figure}[htp]
  \centering
  \includegraphics[scale=0.4]{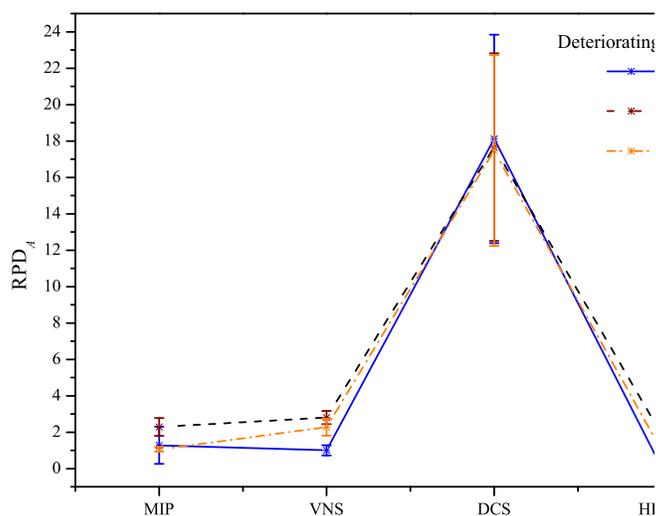}\\
  \caption{Average $\mathrm{RPD}_A$ Means plot and the Tukey HSD intervals (at the 95\% confidence level) for the tested algorithms (small instances)}\label{fig:fig_RPDSI_ANOVA}
\end{figure}

In order to evaluate the algorithms when the size of the problem is increased, experiments are also performed using  large-sized  instances. Since the CPLEX already consumes a lot of computational time in solving small-sized instances, it is excluded in this case. The results of large-sized instances obtained by all the meta-heuristic algorithms are listed in table \ref{tab:tab_ARPDRT_LI}.
In the table, the best RPD value for a specific instance delivered by an algorithm is boldfaced.
We can see that the average results ($\mathrm{RPD}_{A}$)  of the HDCS are {\em always} better than those of the DCS and the VNS.
However, there are 15 instances in which the quality of the best solutions ($\mathrm{RPD}_{B}$) obtained by the HDCS is poorer than that of the VNS. This might be due to the fact that more neighborhood structures are used in the VNS. Nevertheless, {the worst solutions ($\mathrm{RPD}_W$) delivered by the HDCS are significantly better than those of the DCS and the VNS. }

In order to validate the statistical significance of the observed difference in solution quality by different algorithms, we apply again the one-way ANOVA as in the previous scenario. Figure \ref{fig:fig_RPDLI_ANOVA} shows the means plot with the Tukey HSD intervals at 95\% confidence level  for each   algorithm under different deteriorating intervals. It is clearly seen that, the HDCS is statistically better than   the other two algorithms, as  the HDCS consistently and statistically possesses lower Tukey interval compared to the other two algorithms for each of the three deteriorating intervals $\mathrm{H}_1$, $\mathrm{H}_2$ and $\mathrm{H}_3$.

In Table \ref{tab:tab_ARPDRT_LI}, from the row of ``Average", it can be concluded that the average time taken by the DCS algorithm is the shortest, while the average time needed for the HDCS algorithm is the longest and the time consumed by the VNS lies between the two algorithms.
This is due to the fact that, there are more than one elite solution being improved by the  variable neighborhood descent in the HDCS.
{But the computational time of the HDCS is  acceptable for obtaining a good quality solution to the underlying problem.}

Finally, the effect of deteriorating intervals on the performance of the proposed algorithm is analyzed. As shown in Figure \ref{fig:fig_RPDA_LI_HDCS}, the performance of the HDCS somewhat depends on a deteriorating interval, but it is not too much significant. When the deteriorating dates are generated from $\mathrm{H}_3$, the mean RPD values delivered by the HDCS is the smallest.

\begin{table}[htp]\scriptsize
  \centering
  \caption{Computational results of large instances for the algorithms}
    \begin{tabular}{cccccccccccccc}
    \toprule
    \multicolumn{2}{c}{\multirow{2}[4]{*}{Instance}} & \multicolumn{4}{c}{VNS}       & \multicolumn{4}{c}{DCS}       & \multicolumn{4}{c}{HDCS} \\
          \cmidrule[0.05em](r){3-6}
    \cmidrule[0.05em](lr){7-10}
    \cmidrule[0.05em](lr){11-14}
    \multicolumn{2}{c}{} & $\mathrm{RPD}_B$  & $\mathrm{RPD}_A$  & $\mathrm{RPD}_W$  & Time(s) & $\mathrm{RPD}_B$  & $\mathrm{RPD}_A$  & $\mathrm{RPD}_W$  & Time(s) & $\mathrm{RPD}_B$  & $\mathrm{RPD}_A$  & $\mathrm{RPD}_W$  & Time(s) \\
    \midrule
    \multirow{12}[2]{*}{$\mathrm{H}_1$} & 30$\times$4  & 1.49  & 8.50  & 19.09  & 41.63  & 49.61  & 49.63  & 49.63  & 5.16  & \textbf{0.90 } & \textbf{4.65 } & \textbf{8.50 } & 50.38  \\
          & 30$\times$6  & \textbf{0.24 } & 6.98  & 15.99  & 39.69  & 43.71  & 43.99  & 44.06  & 5.27  & 1.57  & \textbf{3.58 } & \textbf{6.22 } & 51.10  \\
          & 30$\times$8  & \textbf{0.89 } & 8.05  & 16.09  & 39.12  & 45.99  & 46.34  & 46.50  & 6.48  & 1.60  & \textbf{4.06 } & \textbf{6.58 } & 46.58  \\
          & 40$\times$4  & \textbf{0.86 } & 10.31  & 24.61  & 95.31  & 65.89  & 67.90  & 68.64  & 7.40  & 1.98  & \textbf{6.09 } & \textbf{10.23 } & 112.20  \\
          & 40$\times$6  & \textbf{1.43 } & 13.75  & 25.76  & 91.41  & 75.04  & 80.10  & 81.40  & 7.92  & 3.79  & \textbf{9.16 } & \textbf{14.80 } & 110.94  \\
          & 40$\times$8  & 2.48  & 10.42  & 19.29  & 78.66  & 50.46  & 51.35  & 51.57  & 7.32  & \textbf{1.06 } & \textbf{5.43 } & \textbf{9.75 } & 102.04  \\
          & 50$\times$4  & 1.99  & 23.64  & 29.63  & 164.51  & 94.65  & 96.52  & 97.01  & 8.76  & \textbf{1.42 } & \textbf{7.40 } & \textbf{15.21 } & 221.95  \\
          & 50$\times$6  & \textbf{3.62 } & 19.39  & 42.70  & 155.02  & 53.94  & 58.40  & 66.08  & 8.68  & 3.75  & \textbf{11.96 } & \textbf{20.18 } & 217.16  \\
          & 50$\times$8  & 3.52  & 16.27  & 30.23  & 151.72  & 76.00  & 76.25  & 76.31  & 8.44  & \textbf{1.79 } & \textbf{9.22 } & \textbf{16.07 } & 183.92  \\
          & 60$\times$4  & 6.85  & 18.67  & 30.55  & 349.18  & 71.67  & 76.47  & 79.71  & 10.48  & \textbf{1.64 } & \textbf{10.92 } & \textbf{18.80 } & 382.93  \\
          & 60$\times$6  & 2.63  & 12.71  & 21.50  & 307.53  & 63.89  & 63.89  & 63.89  & 10.44  & \textbf{2.53 } & \textbf{7.67 } & \textbf{13.38 } & 315.90  \\
          & 60$\times$8  & 3.05  & 11.21  & 19.78  & 316.57  & 57.48  & 57.48  & 57.48  & 10.20  & \textbf{0.60 } & \textbf{6.37 } & \textbf{12.23 } & 337.98  \\
    \multicolumn{2}{c}{Average} & 2.42  & 13.33  & 24.60  & 152.53  & 62.36  & 64.03  & 65.19  & 8.05  & 1.89  & 7.21  & 12.66  & 177.76  \\
          &       &       &       &       &       &       &       &       &       &       &       &       &  \\
    \multirow{12}[0]{*}{$\mathrm{H}_2$} & 30$\times$4  & 1.70  & 8.81  & 16.44  & 43.43  & 35.19  & 35.64  & 35.79  & 5.45  & \textbf{0.96 } & \textbf{5.20 } & \textbf{9.70 } & 46.19  \\
          & 30$\times$6  & 0.94  & 6.12  & 13.54  & 41.03  & 25.84  & 27.36  & 28.15  & 5.44  & \textbf{0.40 } & \textbf{2.51 } & \textbf{4.50 } & 45.30  \\
          & 30$\times$8  & 1.46  & 7.58  & 15.95  & 36.84  & 28.31  & 29.78  & 30.52  & 6.68  & \textbf{0.60 } & \textbf{2.39 } & \textbf{3.88 } & 42.12  \\
          & 40$\times$4  & \textbf{0.81 } & 10.33  & 22.68  & 97.64  & 38.18  & 38.32  & 38.37  & 7.00  & 2.49  & \textbf{7.81 } & \textbf{12.19 } & 104.59  \\
          & 40$\times$6  & \textbf{0.87 } & 10.16  & 23.08  & 90.96  & 37.43  & 38.43  & 39.01  & 6.88  & 2.05  & \textbf{5.67 } & \textbf{10.00 } & 91.51  \\
          & 40$\times$8  & \textbf{0.71 } & 8.66  & 22.33  & 82.68  & 49.37  & 50.00  & 50.15  & 6.64  & 1.66  & \textbf{5.66 } & \textbf{9.81 } & 98.96  \\
          & 50$\times$4  & 2.83  & 17.59  & 38.13  & 156.74  & 60.14  & 61.22  & 61.15  & 9.00  & \textbf{1.34 } & \textbf{7.95 } & \textbf{14.98 } & 184.93  \\
          & 50$\times$6  & \textbf{1.96 } & 11.16  & 21.57  & 168.50  & 26.82  & 26.82  & 26.82  & 8.60  & 1.97  & \textbf{6.39 } & \textbf{11.90 } & 179.63  \\
          & 50$\times$8  & 1.99  & 15.06  & 27.47  & 166.81  & 38.17  & 39.24  & 40.06  & 8.84  & \textbf{1.96 } & \textbf{7.77 } & \textbf{13.31 } & 168.26  \\
          & 60$\times$4  & 3.04  & 10.77  & 18.16  & 310.91  & 32.29  & 32.29  & 32.29  & 10.44  & \textbf{0.97 } & \textbf{6.79 } & \textbf{12.50 } & 316.43  \\
          & 60$\times$6  & 4.55  & 15.73  & 28.03  & 288.22  & 44.30  & 44.30  & 44.30  & 10.36  & \textbf{1.07 } & \textbf{6.62 } & \textbf{11.41 } & 298.84  \\
          & 60$\times$8  & 3.81  & 8.96  & 19.10  & 262.49  & 20.31  & 20.31  & 20.31  & 10.24  & \textbf{1.61 } & \textbf{5.04 } & \textbf{8.68 } & 270.85  \\
    \multicolumn{2}{c}{Average} & 2.06  & 10.91  & 22.21  & 145.52  & 36.36  & 36.98  & 37.24  & 7.96  & 1.42  & 5.82  & 10.24  & 153.97  \\
          &       &       &       &       &       &       &       &       &       &       &       &       &  \\
    \multirow{12}[2]{*}{$\mathrm{H}_3$} & 30$\times$4  & 2.23  & 8.59  & 16.09  & 42.96  & 53.71  & 54.57  & 54.94  & 5.82  & \textbf{0.25 } & \textbf{3.50 } & \textbf{5.76 } & 50.39  \\
          & 30$\times$6  & \textbf{1.27 } & 9.44  & 17.66  & 41.98  & 51.20  & 52.43  & 53.04  & 6.12  & 2.08  & \textbf{4.31 } & \textbf{7.04 } & 46.47  \\
          & 30$\times$8  & 1.69  & 7.24  & 15.42  & 36.66  & 50.18  & 51.67  & 52.44  & 6.16  & \textbf{0.35 } & \textbf{1.86 } & \textbf{3.15 } & 41.30  \\
          & 40$\times$4  & 1.92  & 12.24  & 22.62  & 97.18  & 57.73  & 57.95  & 58.01  & 6.96  & \textbf{0.99 } & \textbf{5.89 } & \textbf{11.00 } & 115.97  \\
          & 40$\times$6  & \textbf{0.92 } & 9.67  & 19.79  & 88.09  & 47.61  & 48.74  & 49.02  & 6.92  & 1.88  & \textbf{5.67 } & \textbf{10.13 } & 103.95  \\
          & 40$\times$8  & \textbf{0.63 } & 8.63  & 18.36  & 84.11  & 66.75  & 67.53  & 67.75  & 6.72  & 1.20  & \textbf{5.62 } & \textbf{10.45 } & 98.52  \\
          & 50$\times$4  & \textbf{3.12 } & 16.66  & 34.57  & 165.67  & 59.92  & 62.83  & 67.96  & 8.68  & 3.66  & \textbf{10.39 } & \textbf{17.57 } & 199.58  \\
          & 50$\times$6  & \textbf{1.01 } & 10.53  & 22.87  & 152.81  & 45.53  & 45.72  & 45.77  & 8.52  & 1.68  & \textbf{6.67 } & \textbf{11.08 } & 189.81  \\
          & 50$\times$8  & \textbf{0.72 } & 13.99  & 29.39  & 148.59  & 81.57  & 82.41  & 82.62  & 8.52  & 3.48  & \textbf{10.43 } & \textbf{18.53 } & 190.57  \\
          & 60$\times$4  & 4.31  & 12.80  & 20.14  & 312.29  & 69.43  & 69.52  & 69.54  & 10.48  & \textbf{1.40 } & \textbf{8.26 } & \textbf{13.26 } & 335.21  \\
          & 60$\times$6  & 4.05  & 15.74  & 28.93  & 310.69  & 62.67  & 62.67  & 62.67  & 10.36  & \textbf{1.83 } & \textbf{7.78 } & \textbf{12.62 } & 332.11  \\
          & 60$\times$8  & 5.32  & 17.55  & 30.32  & 268.94  & 73.55  & 73.55  & 73.55  & 10.20  & \textbf{2.06 } & \textbf{10.88 } & \textbf{18.89 } & 299.49  \\
    \multicolumn{2}{c}{Average} & 2.27  & 11.92  & 23.01  & 145.83  & 59.99  & 60.80  & 61.44  & 7.95  & 1.74  & 6.77  & 11.62  & 166.95  \\
    \bottomrule
    \end{tabular}%
  \label{tab:tab_ARPDRT_LI}%
\end{table}%

\begin{figure}[htp]
  \centering
  \includegraphics[scale=0.4]{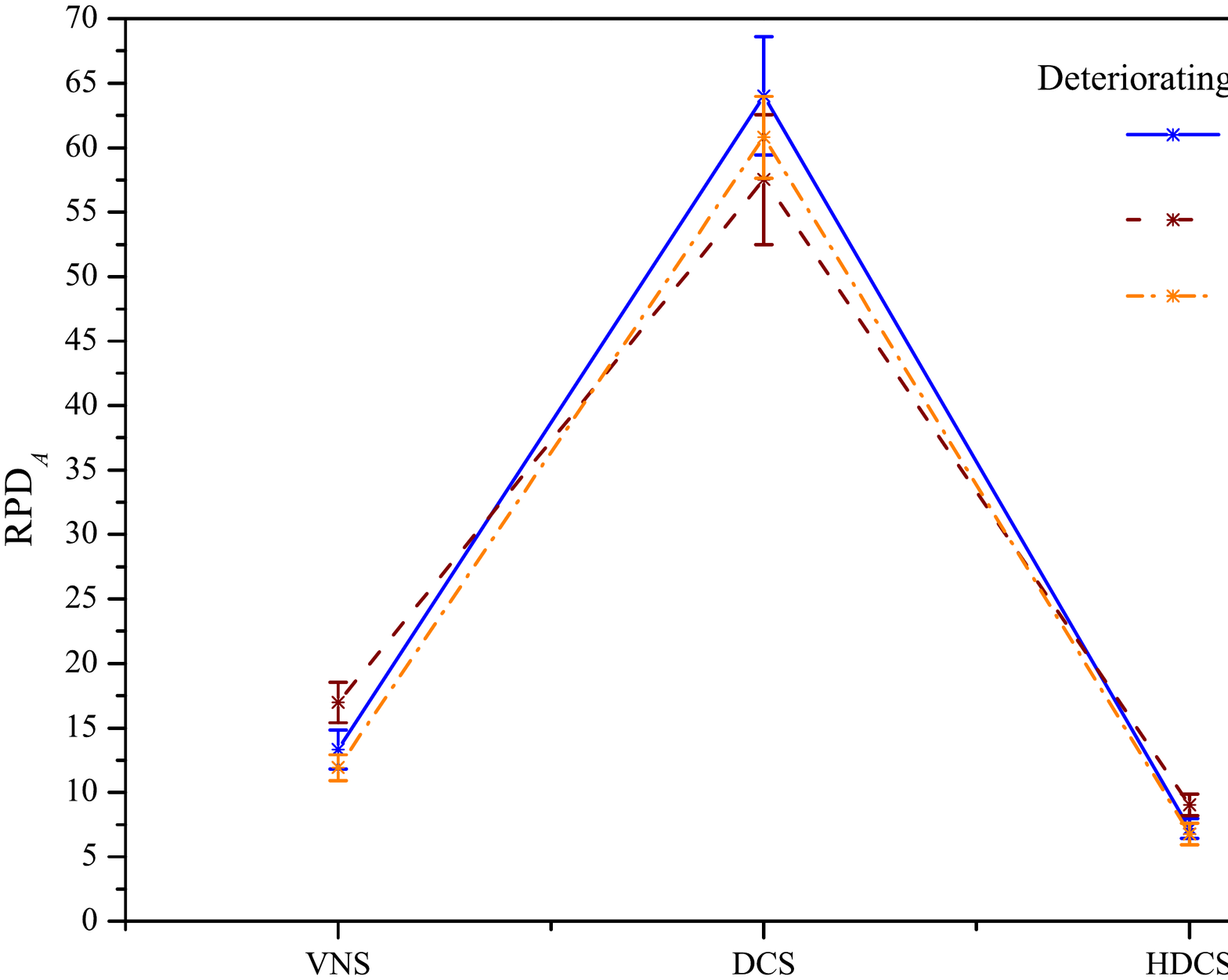}\\
  \caption{Average $\mathrm{RPD}_A$ means plot and the Tukey HSD intervals (at the 95\% confidence level) for the tested algorithms (large instances)}\label{fig:fig_RPDLI_ANOVA}
\end{figure}

\begin{figure}[htp]
  \centering
  \includegraphics[scale=0.4]{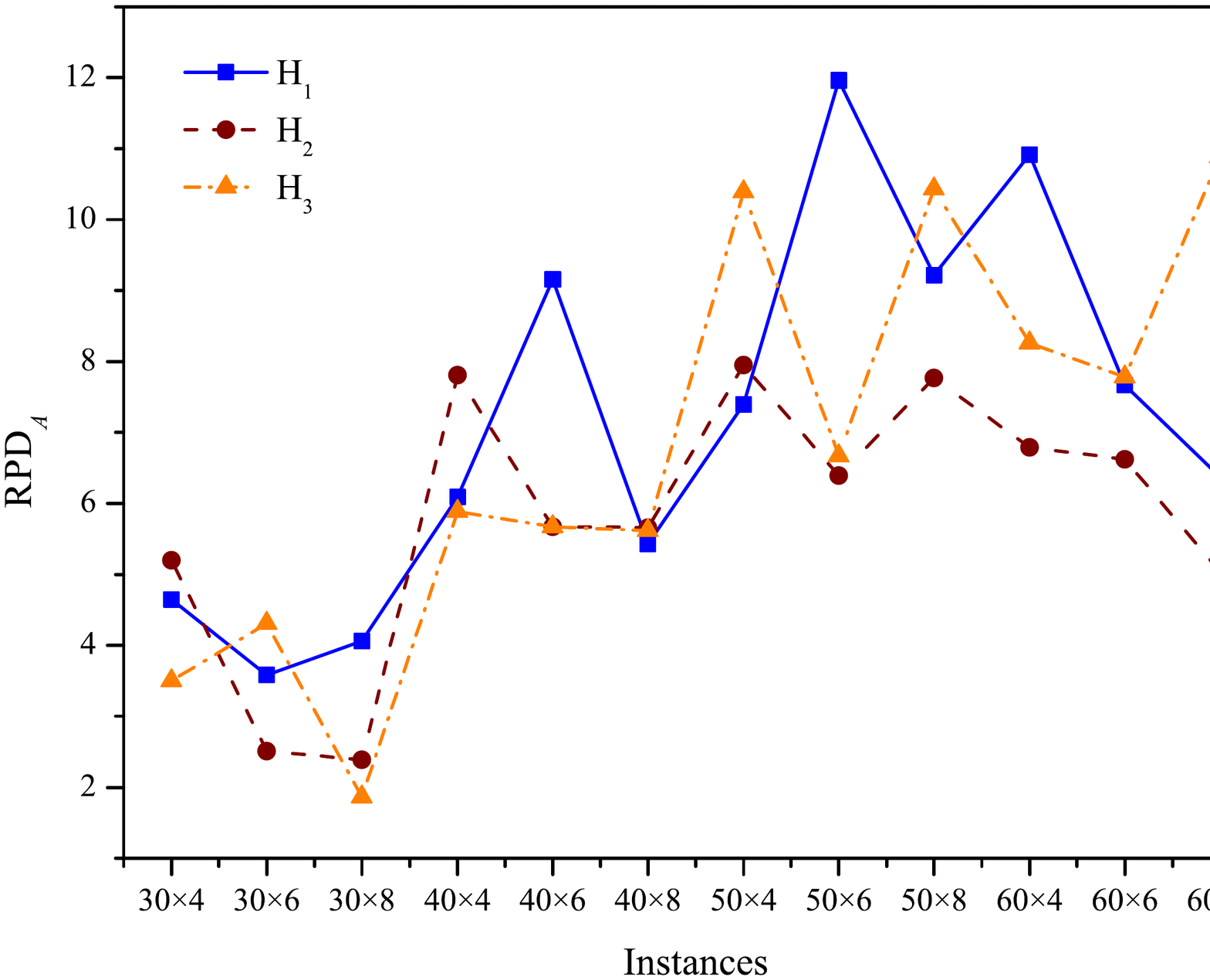}\\
  \caption{The effect of different deteriorating intervals on the performance of the HDCS}\label{fig:fig_RPDA_LI_HDCS}
\end{figure}

\section{Conclusions}\label{sec:sec_conclusion}
In this article, an identical parallel machine scheduling problem with step-deteriorating jobs and sequence-dependent setup times was considered. The objective of this problem is to find a schedule for minimizing the total tardiness. To solve this problem, a mixed-integer programming model was presented to deliver the exact solutions for relatively small-sized   instances. Due to the intractability of the problem, a hybrid discrete cuckoo search algorithm was proposed to obtain a near-optimal solution.
In the proposed algorithm, the operators of the standard CS algorithm are redefined by some discrete operations. To improve the quality of the solution, the variable neighborhood descent  as local search approach is used to refine a set of chosen elite solutions. Moreover, the MBHG heuristic was incorporated into the generation of initial solutions. A set of small to large instances were produced to test the performance of the proposed algorithm.
The computational results show that the hybrid discrete cuckoo search algorithm in general
 provides good average results for test problems,  especially in large-sized instances. Indeed, the HDCS algorithm outperforms the CPLEX solver, the VNS and the DCS algorithms in terms of the average RPD values, shown  in the row `Average' of tables \ref{tab:tab_ARPD_SI} and \ref{tab:tab_ARPDRT_LI}.
 Further studies may concern on applying the proposed algorithm to solve other scheduling problems, such as unrelated parallel machine scheduling problem, the flow shop problem and their stochastic versions. In addition, some multi-objective scheduling problems may be  considered in future work.

\section*{Acknowledgments}

This work was partially supported by the National Natural Science
Foundation of China (No.51175442 and 51205328) and the Fundamental Research Funds
for the Central Universities (No. 2010ZT03; SWJTU09CX022).


\bibliographystyle{gENO}
 \newcommand{\noop}[1]{}

\label{lastpage}

\end{document}